\documentclass []{svmult}
\usepackage[utf8]{inputenc} 
\usepackage[T1]{fontenc} 
\usepackage{amsmath}
\usepackage{amsfonts}
 \usepackage{mathptmx}
 \usepackage{multicol} 
 \usepackage{pstricks}
\usepackage{graphicx}
\usepackage{epstopdf}
\newlength{\oldparindent}
\newcommand{\un}{{\bf 1}}

\newcommand{\cL}{{\mathbb {L}}}%espace L^p

\newcommand{\bpf}{\begin{proof}}
\newcommand{\epf}{ \end{proof} \medskip}

\newcommand{\benum}{\begin{enumerate}}
\newcommand{\eenum}{\end{enumerate}}

\newcommand{\bitem}{\begin{itemize}}
\newcommand{\eitem}{\end{itemize}}

\newcommand{\brmq}{\begin{remark}}
\newcommand{\ermq}{\end{remark}}

\newcommand{\brmqs}{\begin{rmqs}}
\newcommand{\ermqs}{\end{rmqs}}

\newcommand{\bapp}{\begin{application}}
\newcommand{\eapp}{\end{application}}

\newcommand{\bapps}{\begin{applications}}
\newcommand{\eapps}{\end{applications}}

\newcommand{\bdefi}{\begin{definition}}
\newcommand{\edefi}{\end{definition}}

\newcommand{\beq}{\begin{equation}}
\newcommand{\eeq}{\end{equation}}

\def\bpm{\begin{pmatrix}}
\def\epm{\end{pmatrix}}

\newcommand{\bcas}{\begin{cases}}
\newcommand{\ecas}{\end{cases}}

\newcommand{\bex}{\begin{exemp}}
\newcommand{\eex}{\end{exemp}}

\newcommand{\bexs}{\begin{exemps}}
\newcommand{\eexs}{\end{exemps}}

\newcommand{\bprop}{\begin{proposition}}
\newcommand{\eprop}{\end{proposition}}

\newcommand{\bthm}{\begin{theorem}}
\newcommand{\ethm}{\end{theorem}}

\newcommand{\bcor}{\begin{corollary}}
\newcommand{\ecor}{\end{corollary}}

\newcommand{\blem}{\begin{lemma}}
\newcommand{\elem}{\end{lemma}}

\newcommand{\beqna}{\begin{eqnarray}}
\newcommand{\eeqna}{\end{eqnarray}}

\newcommand{\beqnas}{\begin{eqnarray*}}
\newcommand{\eeqnas}{\end{eqnarray*}}

\newcommand{\cA}{{\mathcal A}}

\definecolor{green}{rgb}{0,.7,.2}
\definecolor{orange}{rgb}{0.9,.5,0}

\def\I{1\!{\rm l}}
\def\E{{\mathbb{E}}}%esperance
\def\tr{\textmd{trace}\,}

\def\det{{ \rm{det}}}  %parentheses have been corrected
\def\Id{{\rm{Id}}} %parentheses have been corrected

\def\cA{{\mathcal A }}
\def\cB{{\mathcal B }}
\def\cC{{\mathcal C}}
\def\cD{{\mathcal D}}

\def\cG{{\mathcal  G}}

\def\cJ{{\mathcal  J}}

\def\cL{{\mathcal L }}
\def\cM{{\mathcal  M}}

\def\cP{{\mathcal P }}

\def\cU{{\mathcal  U}}

\def\bbB{{\mathbb{B}}}
\def\bbC{{\mathbb{C}}}%complexes

\def\bbE{{\mathbb{E}}}%esperance

\def\bbL{{\mathbb{L}}}

\newcommand{\bbN}{{\mathbb {N}}}% entiers
\newcommand{\bbR}{{\mathbb {R}}}%reels
\newcommand{\bbS}{{\mathbb {S}}} % sphere 
\def\ccL{{\mathbb{L}}}

\def\discrim{{\rm disc}} 

\author{Dominique Bakry \and  Olfa Zribi}
\institute{D. Bakry \and O. Zribi \at Institut de Mathématiques, Université P. Sabatier, 
118 
route de Narbonne, 31062 Toulouse, FRANCE \email{Dominique.Bakry@math.univ-toulouse.fr , Olfa.Zribi@math.univ-toulouse.fr}}

\makeatletter \renewcommand{\@oddfoot}{\sl \small
 \hfil \thepage\hfil \today}
\renewcommand{\@oddhead}{\sl \small
 \hfil }
 \makeatother

\title{$h$-transforms and orthogonal polynomials
}

\begin{document}
\abstract{ We describe some examples of classical and explicit $h$-transforms as particular cases of a general mechanism, which is related to the existence of symmetric diffusion operators having orthogonal polynomials as spectral decomposition.}



\maketitle
\section{Introduction}

When the first author of this paper was a student, he was attending the DEA course of Marc Yor, about Brownian motions and the many  laws that one would compute explicitly for various transformations on the trajectories. It looked like magic, and  was indeed. In particular, the fact that conditioning a real Brownian motion  to remain positive would turn it into a Bessel process in dimension 3, that is the norm of a 3-dimensional Brownian motion, seemed miraculous.  Of course, there are much more striking identities concerning the laws of Brownian motion that one may find in the numerous papers or books of Marc Yor (see~\cite{MansuyYorBook2008} for a large collection of such examples). The same kind of conditioning appears in  many similar    situations, and specially in geometric models.   This is related to the fact that we then have explicit $h$ (or Doob)- transforms.

This relation between conditioning and $h$-transform was  first put forward by  J.L. Doob~\cite{Doob57}, and is described in full generality in Doob's book~\cite{DoobBook84}. However, this kind of conditioning has been extended in various contexts,  and very reader friendly explained by Marc Yor and his co-authors, in particular in~\cite{RoynetteYor2009, NajnudelRoynetteYor09}. The fact that conditioning a $d$-dimensional model to remain in some set produces a new model in the same family (whatever the meaning of "family"), moreover  with dimension $d+2$,  appears to be a general feature worth to be further understood. It turns out that   the most known models have a common explanation, due to an underlying  structure related to orthogonal polynomials. The scope of this short note is to shed light on  these connections. 

The paper is organized as follows. In Section~\ref{sec.symm.diff}, we present the langage of symmetric diffusion operators that we shall use in the core of the text, and explain what $h$-transforms are.  Section~\ref{sec.examples} gives a few classical and known examples (some of them less well known indeed). They all follow the same scheme,  explained in Section~\ref{sec.orthogonal.poly}, which provides the general framework, related to the study of orthogonal polynomials which are eigenvectors of diffusion operators.  The last Section~\ref{sec.further.examples} provides further examples, as applications of the main result, inspired from random matrix theory.

\section{Symmetric diffusion operators, images and  and $h$-transforms\label{sec.symm.diff}}

\subsection{Symmetric diffusion operators}
We give here a brief account of the tools and notations that we shall be using throughout this paper, most of them following the general setting described in~\cite{bglbook}. A symmetric diffusion process $(\xi_t)$ on a measurable space $E$  may be described by its generator $\cL$, acting on a good algebra  $\cA$  of real valued functions (we shall be more precise about this below). The diffusion property is described through the so-called change of variable formula. Namely, whenever $f= (f_1, \cdots, f_p)\in \cA^p$, and if $\Phi : \bbR^p\mapsto \bbR$ is a smooth function such that $\Phi(f)\in \cA$ together with $\partial_i\Phi(f)$  and $\partial_{ij}\Phi(f)$,  $\forall i,j= 1 \cdots n$,  then  
\beq\label{change.var}\cL(\Phi(f))= \sum_i \partial_i \Phi(f) \cL(f_i) +\sum_{ij} \partial_{ij}\Phi(f) \Gamma(f_i,f_j),\eeq where $\Gamma(f,g)$ is the square field operator (or carré du champ), defined on the algebra $\cA$ through
$$\Gamma(f,g)= \frac{1}{2}\big(\cL(fg)-f\cL(g)-g\cL(f)\big).$$ 

This change of variable formula~\eqref{change.var} is some "abstract" way of describing a second order differential operator with no 0-order term.
It turns out that the operators associated with diffusion processes satisfy
$\Gamma(f,f)\geq 0$ for any $f\in \cA$, and that the operator $\Gamma$ is a first order differential operator in each of its argument, that is, with the same conditions as before,
\beq \label{change.var.Gamma} \Gamma(\Phi(f), g))= \sum_i \partial_i\Phi(f) \Gamma(f_i,g),\eeq

In most cases, our set $E$ is an open subset $\Omega\subset  \bbR^n$, and the algebra $\cA$ is the set of smooth (that is $\cC^\infty$) functions $\Omega\mapsto \bbR$. Then, using formula~\eqref{change.var} for a smooth function $f : \Omega\mapsto \bbR$  instead of $\Phi$ and  $ (x_1, \cdots, x_n)$ instead of  $(f_1, \cdots, f_n)$, we  see that 
$\cL$ may be written as 
\beq\label{diff.coord}\cL(f)= \sum_{ij} g^{ij}(x)\partial^2_{ij} f+ \sum_i b^i(x)\partial_i f,\eeq and   similarly 
$$\Gamma(f,g)= \sum_{ij} g^{ij}(x) \partial_i f\partial_j g.$$ In this system of coordinates, 
$g^{ij}= \Gamma(x_i,x_j)$ and $b^i= \cL(x_i)$.
The positivity of the operator $\Gamma$ just says   that the symmetric matrix $(g^{ij})(x)$ is non negative for any $x\in \Omega$, which is usually translated into the fact that the operator is semi-elliptic. In the same way, the absence of constant term translates into the fact that for the constant function $\un$, that we always assume to belong to the set $\cA$, one has $\cL(\un)=0$, which is an easy consequence of~(\ref{diff.coord}). 

It is not always wise to restrict to  diffusion operators defined on some open subsets of $\bbR^n$. We may have to deal with operators defined on manifolds, in which  case one may describe the same objects in a local system of coordinates. However, using such local system of  coordinates in not  a good idea. In Section~\ref{subsec.matrix.jacobi} for example, we shall consider the group $SO(d)$ of $d$-dimensional orthogonal matrices. The natural algebra $\cA$ that we want to use  is then the algebra of polynomial functions in the entries $(m_{ij})$ of the matrix, and the natural functions $\Phi$ acting on it are the polynomial functions. Since the polynomial structure will play an important rôle in our computations, it is not wise in this context to consider local system of coordinates (the entries of the matrix cannot play this rôle since they are related through algebraic relations).

Coming back to the general situation, the link between the process $(\xi_t)$  and the operator $\cL$ is that, for any $f\in \cA$,  
$$f(\xi_t)- f(\xi_0)-\int_0^t \cL(f)(\xi_s) ds$$ is a local martingale, and this is enough to describe the law of the process starting from some initial point $\xi_0=x\in E$, provided the set of functions $\cA$ is large enough, for example when $\cA$ contains a core of the so-called domain of the operator $\cL$, see  \cite{bglbook}, chapter 3,  for more details.

The law of a single variable $\xi_t$, when $\xi_0=x$, is then described by a Markov operator $P_t$, as 
$$P_t(f)(x)= \bbE_x(f(\xi_t)),$$ and, at least at a formal level, $P_t= \exp(t\cL)$ is the semigroup generated by $\cL$. 

In most of the cases that we are interested in, the operator $\cL$ will be symmetric in some $\ccL^2(\mu)$ space. That is, for some subset $\cA_0$ of $\cA$, which is rich enough to describe $P_t$ from the knowledge of $\cL$ (technically, as mentioned above, a core in the domain $\cD(\cL)$), one has, for $f,g$ in  $\cA_0$
$$\int f\cL(g)\,d\mu= \int g \cL(f)\, d\mu.$$
This translates into the integration by parts formula
\beq \label{eq.ipp} \int f\cL(g)\,d\mu=- \int \Gamma(f,g)\, d\mu.\eeq

For an operator given in an open set $\Omega\subset\bbR^n$ by the formula~(\ref{diff.coord}), and when the coefficients $g^{ij}$ and $b^i$ are smooth, one may identify the density $\rho(x)$ of the measure $\mu$, when $\rho(x)>0$,  by the formula
$$\cL(f) = \frac{1}{\rho(x)} \sum_{ij} \partial_i (\rho g^{ij} \partial_j f),$$
which gives 
\beq\label{eq.density}b^i =\sum_j  (g^{ij}\partial_j \log \rho + \partial_j g^{ij}),\eeq an easy way to recover $\rho$ up to a multiplicative constant  provided $(g^{ij})$ is non degenerate, that is when $\cL$ is elliptic. We call this measure $\mu$  the reversible measure.  Indeed, whenever the measure $\mu$ is a probability measure, and under this symmetry property, then the associated process $(\xi_t)$ has the property that, whenever the law of $\xi_0$ is $\mu$, then for any $t>0$ the law of $(\xi_{t-s}, s\in [0,t])$ is identical to the law of $(\xi_s, s\in [0,t])$. This justifies in this case  the name "reversible", which we keep in  the infinite mass case, following \cite{bglbook}. 

 Through the integration by parts formula, the operator $\cL$ (and therefore the process and the semigroup  themselves,  provided we know something about a core in the domain), is entirely described by the triple $(\Omega, \Gamma, \mu)$, called a Markov triple in~\cite{bglbook}.

Thanks to the change of variable formula~(\ref{change.var}), it is enough to describe an operator in a given system of coordinates $(x^i)$ to describe $\cL(x^i)= b^i$ and $\Gamma(x^i, x^j)= g^{ij}$. Indeed, this determines  $\cL(\Phi(x^i))$, for any $\Phi$ at least $\cC^2$. As outlined earlier, we do not even require that these functions $x^i$  form a coordinate system. They may be redundant (that is more variables than really necessary, as for example in the $SO(d)$ mentioned above), or not sufficient, provided the computed expressions depend only on those variables,  as we do  for example in Section~\ref{sec.further.examples}.

 Moreover, it may be convenient in even dimension to use complex variables, that is,  for a pair $(x,y)$ of functions in the domain, to set $z= x+iy$ and describe $\cL(z)= \cL(x)+ i \cL(y)$, $\Gamma(z,z)= \Gamma(x,x)-\Gamma(y,y)+2i \Gamma(x,y)$ and $\Gamma(z, \bar z)= \Gamma(x,x)+ \Gamma(y,y)$, and similarly for many pairs of real variables, or a pair of a real variable and a complex one. This will be used for example in paragraphs~\ref{subsec.deltoid} and~\ref{sec.Weyl.chamber}. However, we shall be careful in this case  to apply $\cL$  only to polynomial functions in the variables $(x,y)$, replacing $x$ by $\frac{1}{2}(z+ \bar z)$ and $y$ by $\frac{1}{2i} (z- \bar z)$. Then, the various change of variable formulae (on $\cL$ and $\Gamma$) apply when considering $z$ and $\bar z$ as independent variables.

As we already mentioned, it may happen that we can find some functions $X_i, i= 1, \cdots, k$ such that, for any $i$, $\cL(X_i)$ depend only on $(X_1, \cdots, X_k)$ and that the same is true  for $\Gamma(X_i,X_j)$  for any pair $(i,j)$. Then, writing $X= (X_1, \cdots, X_k)\in \bbR^k$,  setting  $B^i(X)= \cL(X^i)$ and $G^{ij}(X) = \Gamma(X_i,X_j)$, one writes for any smooth function $\Phi : \bbR^k\mapsto \bbR$,
$\cL\big(\Phi(X)\big)= \hat \cL (\Phi)(X)$, where 
$$\hat\cL = \sum_{ij} G^{ij}(X) \partial^2_{ij} + \sum_i B_i(X) \partial_i ,$$ which is a direct consequence of formula~(\ref{change.var}). When such happens, the image of the process $(\xi_t)$  with generator $\cL$  under the map $X$ is again a diffusion process $(\hat \xi_t)$ with generator $\hat \cL$. In this situation, we say that $\hat \cL$ is the image of $\cL$ through the map $X$.

Some caution should be taken in this assertion concerning the domains of the operators, but in the examples below all this will be quite clear (our operators will mostly  act on polynomials). When $\cL$ is symmetric with respect to some probability measure $\mu$, then $\hat \cL$ is symmetric with respect to the image measure $\hat \mu$ of $\mu$ through $X$. With the help of formula~(\ref{eq.density}), it may be an efficient way to compute $\hat \mu$.

\subsection{$h$-tranforms}

Given some  diffusion  operator  $\cL$ on some open set in $\bbR^d$, 
 we may  sometimes find  an  explicit function $h$, defined on some subset $\Omega_1$ of $\Omega$,  with values in $(0, \infty)$ such that $\cL(h) = \lambda h$, for some real parameter $\lambda>0$.  We then  look at  the new operator $\cL^{(h)}$, acting on functions defined on $\Omega_1$,  described as 
$$\cL^{(h)}(f)= \frac{1}{h}\cL(hf)-\lambda f$$ is another diffusion operator with the same square field operator than $\cL$. This is the so-called $h$ (or Doob's) transform, see~\cite{Doob57,DoobBook84, bglbook}.
Indeed, thanks to formula~(\ref{change.var}), one has
$$\cL^{(h)}(f) = \cL(f) + 2\Gamma(\log  h, f).$$
 When the operator $\cL$ is symmetric with respect to some measure $\mu$, then $\cL^{(h)}$ is symmetric with respect to $d\mu_h=h^{2}d\mu$.

Considering functions with support in $\Omega_1$, the application $D: f\mapsto hf$ is an isometry between $\ccL^2(\mu_h)$ and $\ccL^2(\mu)$.  It is worth to observe that $\cL^{(h)}= D^{-1} (\cL-\lambda \Id) D$ :  every spectral property (discreteness of the spectrum, cores, etc.) is preserved through this transformation.

 For example, if $f\in \ccL^2(\mu)$ is an eigenvector of $\cL$ with eigenvalue $-\lambda_1$, then $f/h$ is an eigenvector of $\cL^{(h)}$ with eigenvalue $-(\lambda_1+\lambda)$.

Also,  at least formally, for the semigroup $P_t^{(h)}$ associated with $\cL^{(h)}$, one has
$$P_t^{(h)}(f)= e^{-\lambda t} \frac{1}{h} P_t (hf).$$

In general, one looks for positive functions $h$ which vanish at the boundary of $\Omega_1$,  and there is a unique such function $h$  satisfying $\cL(h)= -\lambda h$, usually called the ground state  for $\cL$ on $\Omega_1$. 
This situation appears in general in the following context.
 When $\cL$ is elliptic on $\Omega\subset \bbR^n$, and whenever $\Omega_1$ is bounded, with $\bar \Omega_1 \subset \Omega$, there one may consider the restriction of $\cL$ on $\Omega_1$. If we impose Dirichlet boundary conditions, then the spectrum of this operator consists of a discrete sequence $0>\lambda_0 > \lambda_1 \geq \cdots \geq \lambda_n \cdots$.  The eigenvector $h$ associated with $\lambda_0$ is strictly positive in $\Omega_1$ and vanishes on the boundary $\partial \Omega_1$. This is the  required  ground state $h$ of the operator $\cL$ on $\Omega_1$. 

In  probabilistic terms,   the operator $\cL^{(h)}$ is the generator of the process $(\xi_t)$, conditioned to stay forever in the subset $\Omega_1$.  However, this interpretation is not that easy to check in the general diffusion case. We shall not be concerned here with this  probabilist aspect of this transformation, which is quite well documented in the literature (see \cite{DoobBook84} for a complete account on the subject, and also~\cite{RoynetteYor2009, NajnudelRoynetteYor09} for many examples on conditioning), but rather play around some algebraic aspects of it in concrete examples. However, for the sake of completeness, we shall briefly explain the flavor of this conditioning argument in the simplest example of finite discrete Markov chains, where all the analysis for  justification of the arguments involved may be removed. 

For this, let us consider a finite Markov chain  $(X_n)$ on some finite space $E$, with probability transition matrix $P(x,y)$, $(x,y)\in E^2$, which would play the rôle of $P_1$ in the diffusion context. For simplicity, let us assume that $P(x,y)>0$ for any $(x,y)\in E^2$. Consider now a subset  $A\in E$, and look at the  restriction $P_A$ of the matrix $P$ to $A\times A$. The Perron-Frobenius theorem  asserts that there exists a unique eigenvector $V_0$  for $P_A$, associated with a positive eigenvalue $\mu_0$, which is everywhere positive. This eigenvector $V_0$ corresponds to the ground state $h$ described above in the diffusion context. Then, one may look at the  matrix $Q$ on $A\times A$, defined through $$Q(x,y) = \frac{V_0(y)}{\mu_0V_0(x)} P(x,y),$$ which is a Markov matrix on $A\times A$. This Markov matrix $Q$ plays on $A$  the rôle of $\exp(\cL^{(h)})$
when $h$ is the ground state on $\Omega_1$.

Fix now $n>0$ and  $N>n$. Let $A_N$ be the event $(X_0\in A, \cdots, X_N\in A)$.  For the initial Markov chain $(X_n)$ with transition matrix $P$ and  for $X_0= x\in A$, consider now   the law of 
$(X_0, \cdots, X_n)$ conditioned on $A_N$.  When $F(X_0, \cdots, X_n)= f_0(X_0)\cdots f_n(X_n)$, it is quite easy to check that 
$$
\frac{E(F(X_0, \cdots X_n) \un_{A_N})}{E(1_{A_N})}= \frac{1}{Q^N(1/V_0)(x)}\tilde E\Big(F(X_0, \cdots, X_n) Q^{N-n}(1/V_0)(X_n)\Big),$$
 where $\tilde E$ denotes the expectation for the law of a Markov chain with matrix transition $Q$.
 
 Now, using the irreducibility of the Markov matrix $Q$, one sees that,   when $N$ goes to infinity, both $Q^{N-n}(1/V_0)(X_n)$ and $Q^N(1/V_0)(x)$ converge to $\int \frac{1}{V_0} d\nu$, where $\nu$ is the (unique) invariant  measure for the matrix $Q$. In the limit, we recover the interpretation of the transition matrix transition $Q$ as a matrix of the conditioning of the Markov chain $(X_n)$ to stay forever in $A$.

Coming back to the general case, it is worth to observe that, at least formally, the transformation $\cL\mapsto \cL^{(h)}$ is an involution. Indeed, $\cL^{(h)}(\frac{1}{h})= -\frac{\lambda}{h}$ and 
$(\cL^{(h)})^{(1/h)}= \cL$.
However,  in the usual context of ground states, the interpretation of the associated process as a conditioning is more delicate, since $1/h$ converges to infinity at the boundary of the domain $\Omega_1$.

It is not in general easy to  exhibit explicit examples of such ground states  $h$, but there are many very well known examples in the literature. 
We shall show that in the realm of diffusion processes which are associated to families of orthogonal polynomials, there is a generic argument to provide them, and that this family of examples cover most of the known ones, either directly, either as limiting cases. 

\brmq Observe that, beyond the case where $h$ is a positive eigenvector for $\cL$, one may use the same transformation for any positive function $h$. One may then
look at $\cL^{(h)} (f)= \frac{1}{h} \cL(hf)= \cL(f) + 2\Gamma(\log h, f) + V h$, where 
$V= \frac{\cL(h)}{h}$. 
In particular,  with operators in $\bbR^n$ of the form $\cL(f) = \Delta (f) + \nabla \log V\cdot\nabla f$, which have reversible measure $V dx$, one may use $h=  V^{-1/2}$, which transforms in an isospectral way $\cL$ into a Shr\"odinger type operator $\Delta f + Vf$,  associated with Feynman-Kac semigroups. This allows to  remove a gradient vector field, the price to pay is  that one adds  a potential term. This technique may be used to analyse spectral properties of such symmetric diffusion operators through the techniques  used for  Shr\"odinger operator (see~\cite{bglbook}, for example).

\ermq

\section{Some examples\label{sec.examples}}

\subsection{Bessel operators} We start form the Brownian motion in $\bbR$. The operator $\cL$ is given by $\cL(f)= \frac{1}{2} f''$.   Here, $\Gamma(f,f)= \frac{1}{2}f'^2$ and $\mu$ is the Lebesgue measure. If we  consider  $\Omega= (0, \infty)$ and $h= x$, one has $\lambda= 0$ and $\cL^{(h)}(f'')= \frac{1}{2}  (f'' + \frac{2}{x} f')$.  This last operator is a Bessel operator $\cB_3$. More generally, a Bessel process $Bes(n)$ with parameter  $n$ has a  generator in $(0,\infty)$ given by
$$\cB_n(f)= \frac{1}{2} (f'' + \frac{n-1}{x} f'),$$ and it is easily seen, when $n\geq 1$ is an integer,  to be the generator  of $\|B_t\|$, where $(B_t)$ is an $n$-dimensional Brownian motion (indeed, $\cB_n$ is the image of the Laplace operator $\frac{1}{2}\Delta$ under $x\mapsto \|x\|$, in the sense described in  Section~\ref{sec.symm.diff}). This $\cB_3$ operator is also the generator of a real Brownian motion conditioned to remain positive. Observe however that the function $h$ is this case does not vanish at the infinite boundary of the set $(0, \infty)$, and that the probabilistic  interpretation would require some further analysis than the one sketched in the previous section.

From formula~(\ref{eq.density}), it is quite clear that a reversible  measure for the operator $\cB_n$ is $x^{n-1} dx$ on $(0, \infty)$, which for $n\in \bbN^*$, is also, up to a constant,  the image of the Lebesgue measure in $\bbR^n$ through the map $x\mapsto \|x\|$.

This $h$-transform may be extended to the general Bessel operator. Indeed, for any $n>0$, one may consider the function $h_n(x)= x^{2-n}$,  for which $\cB_n(h_n)=0$, and then 
$\cB_n^{(h_n)}= \cB_{4-n}$.

The change of $\cB_n$ into $\cB_{4-n}$ is perhaps more clear if we consider the generator through the change of variable $x\mapsto x^2$, that is if we consider the generator of the process $(\xi_t^2)$ instead of the process $(\xi_t)$ with generator $\cB_n$. A simple change of variable provides the image operator
\beq\label{gen.bess.poly}\hat \cB_n(f)= 2xf'' + nf',\eeq for which the reversible measure has density $\rho(x)= x^{(n-2)/2}$,  and the function $h$ is nothing else than  $1/\rho$.

 Under this form, we shall see that is a particular case of a phenomenon related to orthogonal polynomials, developed in Section~\ref{sec.orthogonal.poly}, although here there are no polynomials involved here, the reversible measure  being infinite.

\brmq It is not hard to observe that for $0<n<2$, the process $(\xi_t)$ with associated generator $\cB_n$, and starting from $x>0$  reaches $0$ in finite time. Then, $\cB_{4-n}$ is the generator of this process conditioned to never reach $0$.  However, it is well known that the Bessel operator is essentially self-adjoint on $(0, \infty)$ as soon as $n>3$ (see \cite{bglbook}, page 98, for example).  This means that the set of smooth function compactly supported in $(0,\infty)$ is dense in the $\ccL^2$ domain of $\cB_n$. Since this is a spectral property, it is preserved through $h$-transform and this also shows that it is  also essentially self adjoint for any $n<1$. In particular, there is a unique symmetric semi-group for  which the generator coincides with $\cB_n$ on the set of smooth compactly supported functions. 
On the other hand, for $1\leq n <2$, since the associated operator hits the boundary in finite time, there are at least two such semigroups with $\cB_n$ as generator acting on smooth functions, compactly supported in $(0, \infty)$ : the one corresponding to the Dirichlet boundary condition, corresponding to the process killed at the boundary $\{x=0\}$, and the one corresponding to the Neuman boundary condition, corresponding to the process reflected  at the boundary. Through $h$-transforms, one sees then that there are also at least two positivity preserving semi groups in the case $2<n\leq 3$, which may be a bit surprising since then the associated process does not touch the boundary. However, although the   Dirichlet semigroup is Markov ($P_t (\un) <\un$), its $h$-transform is Markov ($P_t (\un) =\un$), while the $h$-transform of the Neuman  semigroup (which is Markov), satisfies $P_t(\un)\geq \un$.

\ermq

\subsection{Jacobi operators\label{subsec.jacobi}}

This is perhaps the most celebrated case of known explicit $h$-transform, since it is closely related in some special case to the Fourier transform on an interval. The Jacobi operator on the interval $(-1,1)$ has generator
$$\cJ_{\alpha, \beta}(f)= (1-x^2) f''- \big((\alpha+\beta)x+ \alpha-\beta\big)f'$$ and is symmetric with respect to the Beta  measure on $(-1,1)$ which is 
$C_{\alpha, \beta} (1-x)^{\alpha-1}(1+x)^{\beta-1} dx$, $C_{\alpha, \beta}$ being the normalizing constant. We always assume that $\alpha, \beta >0$. There is a duality through $h$-transforms exchanging $\cJ_{\alpha, \beta}$ and $\cJ_{2-\alpha, 2-\beta}$, the function $h$ being $ (1-x)^{1-\alpha}(1-x)^{1-\beta}$, that is,  as in the Bessel case in the appropriate coordinate system, the inverse of the density measure.

In a similar way that the Bessel process may be described as a norm of a Brownian motion, one may see the symmetric Jacobi  operator  ($\alpha= \beta$) as an image of a spherical Brownian motion in dimension $2\alpha$. Namely, if one considers the unit sphere $\bbS^n$ in $\bbR^{n+1}$, and looks at  the Brownian motion on it (with generator $\Delta_{\bbS^n}$ being the Laplace operator on the sphere), and then one looks at its first component, one gets a 
 process on $(-1,1)$ with generator $\cL^{n/2,n/2}$. (We refer to paragraph~\ref{subsec.spheres} for details about the spherical Laplacian, from which this remark follows easily, see also~\cite{ bglbook,SteinWeiss}). One may also provide a similar description in the asymmetric case, when the parameters $\alpha$ and $\beta$ are half integers. In this case, $\cL_{\alpha, \beta}$ is, up to a factor 4,  the image of the spherical Laplace operator acting on the unit sphere $\bbS^{2\alpha+ 2\beta-1}$ through the function $X : \bbS^{2\alpha+ 2\beta-1}\mapsto [-1,1]$ defined, for $x= (x_1, \cdots, x_{2\alpha+ 2\beta}) \in \bbR^{2\alpha+ 2\beta}$ as
 $$X(x) = -1+ 2\sum_{i=1}^{2\alpha}x_i^2.$$

The operator $\cJ_{\alpha, \beta}$ may be diagonalized in a basis of orthogonal polynomials, namely the Jacobi polynomials.  They are deeply related to the analysis on the Euclidean case in the geometric cases described above. For example, when $\alpha= \beta$ is an half-integer, then,   for each degree $k$, and up to a multiplicative constant,  there exists a unique function on the sphere  which depends only on the first coordinate and which is the restriction to the sphere of an homogeneous  degree $k$ harmonic polynomial in the corresponding Euclidean space : this is the corresponding  degree $k$ Jacobi polynomial (see~\cite{SteinWeiss, bglbook} for more details). In other words, if $P_k(x)$ is one of these Jacobi polynomials with degree $k$ corresponding to the case $\alpha= \beta= n/2$, then the function $(x_1, \cdots, x_{n+1})\mapsto\|x\|^k P_k(\frac{x_1}{\|x\|} )$ is an homogeneous  harmonic polynomial in $\bbR^{n+1}$.  A similar interpretation is valid in the asymmetric case, whenever the parameters $\alpha$ and $\beta $ are half-integers, if one reminds that the eigenvectors of the Laplace operator on the sphere are restriction to the sphere of harmonic homogeneous polynomials in the ambient Euclidean space (see~\cite{SteinWeiss}).

For $\alpha= \beta= 1/2$, $\cJ_{\alpha, \beta} $ this is just the image of the usual operator $f''$ on $(0, \pi)$ through the change of variables $\theta\mapsto \cos(\theta)=x$. More generally, in the variable $\theta$, $\cJ_{\alpha, \beta} $ may be written as
$$\cJ_{\alpha, \beta} = \frac{d^2}{d\theta^2} + \frac{(\alpha+ \beta-1)\cos(\theta) + \alpha-\beta}{\sin(\theta)} \frac{d}{d\theta}.$$

For $\alpha= \beta= 1/2$, corresponding to the arcsine law, the associated orthogonal polynomials $P_n^{1/2,1/2}$ are the Chebyshev polynomials  of the first kind, satisfying
$$P_n^{1/2,1/2}(\cos(\theta))= \cos(n\theta).$$ 
For $\alpha= \beta= 3/2$, corresponding to the semicircle law,  they correspond to the Chebyshev polynomials of the second kind, satisfying the  formula
$$\sin(\theta)P_n^{3/2,3/2}(\cos(\theta))= \sin(n\theta).$$
These formulae indeed reflect the $h$-transform between $\cJ^{1/2,1/2}$ and $\cJ^{3/2,3/2}$. While 
$P_n^{1/2, 1/2}(\cos(\theta))$ is a basis of $\ccL^2\big((0,\pi), dx\big)$ with Neuman boundary conditions, 
$ \sin(\theta)P_n^{3/2,3/2}(\cos(\theta))$ is another basis of  $\ccL^2\big((0,\pi), dx\big)$, corresponding to the Dirichlet boundary condition. This is the image of the eigenvector basis for $\cL^{3/2, 3/2}$ through the inverse $h$ transform,  the function $h$ being in this system of coordinates nothing else than $(\sin\theta)^{-1}$. 

 For $n=1$, one gets the projection of the Brownian motion on the circle, which is locally a Brownian motion on the real line, up to a change of variables. The first coordinate $x_1$  on the sphere plays the rôle of a distance to the point $(1, 0, \cdots, 0)$ (more precisely, $\arccos(x_1)$ is the Riemannian distance on the sphere from $(1, 0, \cdots, 0)$ to any point with first coordinate $x_1$), and we have a complete analogue of the case of the one dimensional Brownian motion.  Namely,
 
 \bprop The Brownian motion on the half interval (identified with the circle) conditioned to never reach the boundaries is, up to a change of variable,  the radial part of a Brownian motion on a 3 dimensional sphere.
 \eprop

\subsection{ Laguerre operators\label{subsec.laguerre}}

This is the family of operator on $(0, \infty)$ with generator 
$$\cL_{(\alpha)} (f) = xf''+ (\alpha-x)f',$$ which is symmetric with respect to the gamma measure   
$$d\mu^{(\alpha)}= C_\alpha x^{\alpha-1} e^{-x} dx.$$
For $\alpha>0$, the Laguerre family of operators is another instance of diffusion operators on the real line which may be diagonalized in a basis of orthogonal polynomials : these polynomials are the Laguerre polynomials, and  are one of the three families, together with Jacobi polynomials and Hermite polynomials, of orthogonal polynomials in dimension 1 which are at the same time eigenvectors of a diffusion operator, see~\cite{BakryMazet03}. The Laguerre operator   is closely related to the Ornstein-Uhlenbeck operator defined in~\eqref{OU-Rn}, and plays for this operator the same rôle that the one  played by Bessel operators for the Euclidean Brownian motion. 

It is indeed quite close to the Bessel generator under the form~\eqref{gen.bess.poly}, and  in fact the Bessel operator may be seen as a limit of Laguerre operators under proper rescaling.  It is also a limit of asymmetric Jacobi operators, also  under proper rescaling (see~\cite{bglbook}).  The function $h= x^{1-\alpha}$ satisfies $\cL_{(\alpha)}(h)= (\alpha-1)h$, and the $h$-tranform of $\cL_{(\alpha)}$ is $\cL_{(2-\alpha)}$. 

 As mentioned above, when  $\alpha$ is a half-integer $n/2$, the Laguerre operator may be seen as the radial part of the Ornstein-Uhlenbeck operator in $\bbR^n$ with generator 
\beq\label{OU-Rn}\cL^{OU}= \Delta -x \nabla,\eeq
 which is symmetric with respect to the standard Gaussian measure.   More precisely, for $\alpha= n/2$, $\cL^{OU}f(\frac{\|x\|^2}{2})= 2\Big(\cL_{(\alpha)}f\Big)(\frac{\|x\|^2}{2})$. It is therefore  an image of the $n$-dimensional  Ornstein-Uhlenbeck  operator in the sense of Section~\ref{sec.symm.diff}.
In other words, the Laguerre process with generator $2\cL_{(n/2)}$ is nothing else than  the squared norm of an Ornstein-Uhlenbeck process in $\bbR^n$. For $\alpha= 1/2$, this corresponds to the modulus of a one dimensional Ornstein-Uhlenbeck,  that is the one dimensional Ornstein-Uhlenbeck operator itself on $(0, \infty)$, and  we get, as the particular case for $n= 1/2$, 
\bprop The law of an Ornstein-Uhlenbeck operator in dimension 1,  conditioned to remain positive is the same as the law of the norm of a 3-dimensional Orntein-Uhlenbeck operator.

\eprop

\subsection{ An example in $\bbR^2$\label{subsec.deltoid}}

The following example, less well known,  had been pointed out by T. Koornwinder \cite{Koorn1}, not exactly under this form of $h$-transform, but in terms of duality between two families of orthogonal polynomials in dimension 2. It shows that the law of a  Brownian motion in the plane, conditioned not to reach the boundaries of an equilateral triangle, has the law of the spectrum of an  Brownian $SU(3)$ matrix. 

This example, closely related to root systems  and reflection groups in the plane, consists in observing the image of a planar Brownian motion reflected along the edges of an equilateral triangle.  This triangle generates a triangular lattice in the plane,  and this image is observed through some function $Z : \bbR^2\mapsto \bbR^2$ which has the property that any function $\bbR^2\mapsto \bbR$ which is invariant under the symmetries among the lines of the lattice is a function of $Z$.  This image of $\bbR^2$ through the function $Z$ is a bounded domain in $\bbR^2$, with boundary the Steiner's hypocycloid.

The Steiner hypocycloid (also called deltoid curve) is the curve obtained in the plane by rotating (from inside) a circle with radius 1 on a circle with radius 3. Is is the boundary of a bounded open  region  in the plane which we call the deltoid domain $\Omega_D$. It is an algebraic curve of degree 4. It's equation may be written in complex coordinates as $\{D(Z, \bar Z)=0\}$, where $D$ is defined in Proposition~\ref{prop.deltoid}.

\begin{figure}[ht]
 \centering		\includegraphics[width=.5\linewidth]{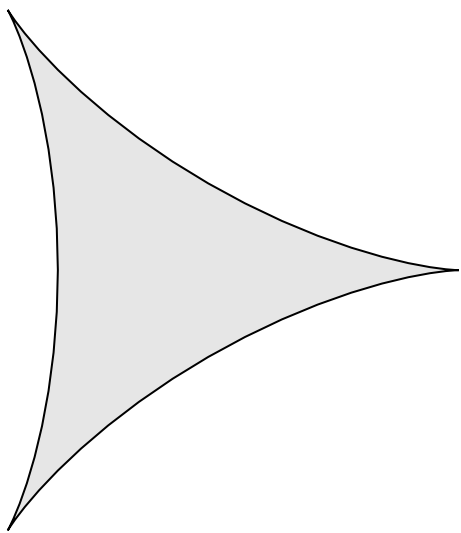}
		\caption{The deltoid domain.}
		\label{fig:Deltoide}
\end{figure}

Consider the following application
$\bbR^2\mapsto \bbR^2$, which is defined as follows. Let $(1, j, \bar j)$ be the  three third  roots of units in the complex plane $\bbC$, and, identifying $\bbR^2$ with $\bbC$, let $Z(z): \bbR^2\mapsto \bbR^2$ be the function 
 $$ Z(z) = \frac{1}{3}\Big(\exp(i (1\cdot z))+ \exp(i(j.\cdot z))+ \exp(i(\bar j\cdot z))\Big),$$ where $z_1\cdot z_2$ denotes the scalar product in $\bbR^2$.
 
 We have
 \bprop ~ Let $L$ be the lattice generated in the plane by the points with coordinates  $M_1=(0, 4 \pi/3)$ and $M_2=(2\pi/3, 2\pi/\sqrt{3})$, and $T$ the (equilateral) triangle with edges $\{(0,0),  M_1,M_2\}$. 
 \benum
 \item The image of $\bbR^2$ under the function $Z$ is the closure $ \bar \Omega_D$ of the deltoid domain.
 
 \item $Z : \bbR^2\mapsto \bbR^2$ is invariant under the symmetries along all the lines of the lattice $L$. Moreover, it is injective on the triangle  $T$.
 \eenum
 
  \eprop
  
  We shall not give a proof of this, which may be checked directly. We refer to~\cite{Zribi2013} for details.  As a consequence, any measurable function $\bbR^2\mapsto \bbR$ which is invariant under the symmetries of $L$ may be written $f(Z)$, for some measurable function  $f : \Omega_D\mapsto \bbR$.
 
The particular choice of this function $Z$ is due to the fact that the Laplace operator in $\bbR^2$ has a nice expression through it. Using complex coordinates as described in Section~\ref{sec.symm.diff},  one has

\bprop \label{eq.deltoid} For the Laplace operator $\Delta$ in $\bbR^2$ and its associated square field operator $\Gamma$, one has
 
 \beq\label{eq2.deltoid}\bcas
 \Gamma(Z,Z) = \bar Z-Z^2, \Gamma(\bar Z, \bar Z)= Z-\bar Z^2, \\ 
 \Gamma(\bar Z, Z)= 1/2(1-Z\bar Z),\\
 \Delta(Z)= -Z, \Delta(\bar Z)= - \bar Z,
 \ecas
\eeq

\eprop

This may be checked directly. One sees that the Laplace operator in $\bbR^2$ has an image through $Z$ in the sense described in Section~\ref{sec.symm.diff}, given  in Proposition~\ref{eq.deltoid}.  This describes the generator of the Brownian motion in the plane,  reflected along the lines of this lattice, coded through this change of variables. One may express the image measure of the Lebesgue measure on the triangle in this system of coordinates. With the help of formula~\eqref{eq.density}, we get 
\bprop \label{prop.deltoid}
Let $D(Z, \bar Z)= \Gamma(Z, \bar Z)^2-\Gamma(Z,Z)\Gamma(\bar Z, \bar Z)$, where $\Gamma$ is given by equation~\eqref{eq2.deltoid}. Then,

\benum 
\item  $D(Z, \bar Z)$ is positive on $\Omega_D$. 

\item $\{D(Z, \bar Z)= 0\}$ is the deltoid curve (that is the boundary of $\Omega_D$).

\item The reversible measure for the image operator described by~\eqref{eq2.deltoid} has density $D(Z,\bar Z)^{-1/2}$ with respect to the Lebesgue measure.

\item  If we write $z_1= \exp(i (1\cdot z))$, $z_2= \exp(i(j.\cdot z))$,  $z_3=  \exp(i(\bar j\cdot z))$, then $$D(Z, \bar Z)= -(z_1-z_2)^2(z_2-z_3)^2(z_3-z_1)^2/(2^23^3).$$

\eenum

\eprop

\brmq \label{rmq.pol.caract.} Observe that thanks to the fact that $|z_i|=1$ and $z_1z_2z_3=1$, the expression $(z_1-z_2)^2(z_2-z_3)^2(z_3-z_1)^2$ is always non positive. Moreover, given  a complex number $Z$ in the deltoid domain $\Omega_D$, there exist three different complex numbers  $(z_1,z_2,z_3)$ with $|z_i|= 1$ and $z_1z_2z_3=1$   such that $Z= \frac{1}{3}(z_1+z_2+z_3)$. They are  unique up to permutation, and are   the solutions of $X^3-3Z X^2+ 3\bar Z  X-1=0$.   Indeed,  for such numbers $z_1,z_2,z_3$, 
$$3\bar Z = \bar z_1+ \bar z_2+ \bar z_3= \frac{1}{z_1} + \frac{1}{z_2}+\frac{1}{z_3}= z_2z_3 +z_1z_3+z_1z_2.$$
\ermq

One may now consider the family of operator $\cL^{(\lambda)}$  defined through
 \beq\label{eq.deltoid.lambda}
 \bcas 
 \Gamma(Z,Z) = \bar Z-Z^2, \Gamma(\bar Z, \bar Z)= Z-\bar Z^2, \\ \Gamma(\bar Z, Z)= 1/2(1-Z\bar Z),\\
\cL^{(\lambda)}(Z)= -\lambda Z, \cL^{(\lambda)}(\bar Z)= -\lambda \bar Z,
\ecas
\eeq 
which is symmetric with respect to the measure $\mu_\lambda= D(Z, \bar Z)^{(2\lambda-5)/6}dZ$, with support the set $\{D(Z,\bar Z) \geq 0\}$ (where $dZ$ is a short hand for the Lebesgue measure in the complex plane) as a direct (although a bit  tedious) computation shows from a direct application of   formula~(\ref{eq.density}) (see Section~\ref{sec.orthogonal.poly} for a proof in a general context which applies in particular here). 

This family of operators plays a rôle similar in this context to the one played by the family $\cJ_{\alpha, \beta}$ for Jacobi polynomials introduced in Section~\ref{subsec.jacobi} or for the family $\cL_{(\alpha)}$ introduced in Section~\ref{subsec.laguerre} for Laguerre polynomials.

This density equation~(\ref{eq.density}) indicates that, for any pair of smooth functions compactly supported in $\{D(Z, \bar Z )>0\}$, the integration by parts~(\ref{eq.ipp}) holds true. Indeed, we have a much stronger result, which extends this formula to any pair of smooth functions defined in a neighborhood of $\bar \Omega$. This relies of some miraculous property of $\partial \Omega$ itself, which has as boundary equation $\{D(Z, \bar Z)=0\}$ and for which
\beq\label{eq.bord.deltoide}\bcas \Gamma(Z, Z)\partial_Z D+ \Gamma(Z, \bar Z)\partial_{\bar Z} D=-3 ZD,\\
 \Gamma(\bar Z , Z)\partial_Z D+ \Gamma(\bar Z, \bar Z)\partial_{\bar Z} D=-3 \bar Z D.
 \ecas
 \eeq
 In particular, $\Gamma(Z, D)$ and $\Gamma(\bar Z, D)$ vanish on $\{D=0\}$. This is a sufficient (and indeed necessary) for the integration by parts~ formula~\eqref{eq.ipp} to be valid  for any pair smooth functions restricted on the set  $\{D\geq 0\}$, in particular for any pair of polynomials (see~\cite{BOZ2013}). Since on the other hand the operator $\cL^{(\lambda)}$ maps polynomials in $(Z,\bar Z)$ into polynomials, without increasing their total degrees, the restriction of $\cL^{(\lambda)}$ on the finite dimensional space of polynomials with total degree less than $k$ is a symmetric operator (with respect to the $\ccL^2(\mu_\lambda)$-Euclidean structure)   on this linear space.  We may therefore find an orthonormal  basis of such polynomials which are eigenvectors for $\cL^{(\lambda)}$, and therefore construct a full orthonormal basis of polynomials made of eigenvectors for $\cL^{(\lambda)}$. 
 
 These polynomials are an example of Jack's polynomials associated with root systems (here the root system $A_2$), see~\cite{Macdo1, DX}, generalized by MacDonald~\cite{Macdo, Macdo1,McDon03}, see also~\cite{KirillovEtingof942,KirillovEtingof95,Cheredn95}, and for which the associated generators are Dunkl operators of various kinds, see~\cite{Heck1,Heck2,Dunkl, Rosler,Sghaier12}. 

For $\lambda= 4$, it turns out that this operator is, up to a scaling factor $8/3$,  the image of the Laplace (Casimir) operator on $SU(3)$ acting on the trace of the matrix. More precisely, on the compact semi-simple Lie group $SU(3)$, we associate to each element  $E$ in the Lie algebra $\cG$ a (right) vector field $X_E$ as follows 
$$X_E(f)(g)= \partial_t (f(ge^{tE}))_{\mid t=0}.$$ Then, one   choses in the Lie algebra  $\cG$ an orthonormal basis $E_i$ for the Killing form (which is negative definite), and we consider the operator $\cL= \sum_i X_{E_i}^2$. This is the canonical Laplace operator on the Lie group, and it commutes with the group action, from left and right : if $L_g(f)(x)= f(xg)$, and $R_g(f)(x)= f(gx)$,  then $\cL L_g= L_g \cL$ and $\cL R_g= R_g\cL$.  For the Casimir operator acting on the entries $(z_{ij})$ of an $SU(d)$ matrix, one  may compute explicitly this operator, and obtain, up to a factor 2,     the following   formulae 
\beq \label{casimir.SU(d)}\bcas \cL^{SU(d)}(z_{kl})= -2\frac{(d-1)(d+1)}{n} z_{kl}, ~\cL^{SU(d)}(\bar z_{kl})= -2\frac{(d-1)(d+1)}{n} \bar z_{kl}\\
\\
 \Gamma^{SU(d)}(z_{kl},z_{rq})= -2z_{kq}z_{rl}+\frac{2}{d}z_{kl}z_{rq}, ~\Gamma(z_{kl},\bar z_{rq})= 2(\delta_{kr}\delta_{lq} -\frac{1}{d}z_{kl}\bar z_{rq}).
 \ecas
 \eeq

 A Brownian motion on $SU(d)$ is a diffusion process which has this Casimir operator as generator   (there are of course many other equivalent definitions of this Brownian motion). 
 
 On $SU(3)$, if one considers the function  $SU(3)\mapsto \bbC$ which to $g\in SU(3)$ associates $Z(g)= \frac{1}{3} \tr(g)$, then one gets for this function $Z$ and for this Casimir operator,  an image operator which is the operator $\frac{8}{3}\cL^{(4)}$, where $\cL^{(\lambda)}$ is defined through equation~(\ref{eq.deltoid.lambda}).  Of course, one may perform the computation directly, or use the method described in paragraph~\ref{sec.Weyl.chamber} to compute  from the operator given of $SU(d)$ through formulas~(\ref{casimir.SU(d)}),  the actions of the generator and the  carré du champ on the characteristic polynomial $P(X)= \det(X\Id-g)$ (see also~\cite{BourgadeHughesNikeYor08} for another approach, together with~\cite{BourgadeYor2009} for nice connections with the Riemann-Zeta function). 

It is worth to observe that functions on $SU(3)$ which depend only on this renormalized trace $Z$ are nothing else but spectral functions. Indeed, if a matrix $g\in SU(3)$ have eigenvalues $(\lambda_1, \lambda_2, \lambda_3)$, with $|\lambda_i|=1$ and $\lambda_1\lambda_2\lambda_3=1$, then a spectral function, that is a symmetric function of $(\lambda_1, \lambda_2, \lambda_3)$, depends only on $\lambda_1+ \lambda_2+ \lambda_3= 3Z$ and,  as observed in Remark~\ref{rmq.pol.caract.}, 
$\lambda_1\lambda_2+ \lambda_2\lambda_3+ \lambda_3\lambda_1= 3\bar Z$.

 Then,  using the function  $D$ which is the determinant of the metric involved in equation~(\ref{eq.bord.deltoide}), one may check directly that 
$$\cL^{(\lambda)}\big(D(Z, \bar Z)^{(5-2\lambda)/6}\big)= (2\lambda-5) D(Z, \bar Z)^{(5-2\lambda)/6},$$ so that one may use the function $h=D(Z, \bar Z)^{(5-2\lambda)/6}$ to perform an $h$ transform on $\cL^{(\lambda)}$ and we obtain 
$$(\cL^{(\lambda)})^{(h)}= \cL^{(5-\lambda)}.$$

Indeed, as we shall see in Section~\ref{sec.orthogonal.poly}, this $h$-transform identity relies only on equation~(\ref{eq.bord.deltoide}). In particular, moving back to the triangle through the inverse function $Z^{-1}$,  for $\lambda=1$, which corresponds to the Brownian motion reflected at the boundaries of the triangular lattice, the $h$ transform is $\cL^{(4)}$, which corresponds to the spectral measure on $SU(3)$.  Then, for this particular case $\lambda=1$, we get
\bprop
A Brownian motion in the equilateral triangle $T$, conditioned to never reach the boundary of the triangle,   has the law of the image under $Z^{-1}$ of the spectrum   of an $SU(3)$ Brownian matrix

\eprop

\subsection{An example in the unit ball in $\bbR^d$\label{subsec.spheres}}

Another example comes from the spherical Brownian motion on the unit sphere  $$\bbS^d= \{(x_1, \cdots, x_{d+1})\in \bbR^{d+1}, \sum_i x_i^2=1\}.$$ To describe the Brownian motion on $\bbS^d$, we look at  its generator, that is  this  the spherical Laplace operator may. It may be described through its action on  the restriction to the sphere of the coordinates $x_i$, seen as   functions $\bbS^{d}\mapsto \bbR$. Then,for the Laplace operator $\Delta^{\bbS^{d}}$ and its associated carré du champ operator $\Gamma$, one has
\beq\label{lapl.sphere}\Delta^{\bbS^{d}}(x_i)= -dx_i, \quad\Gamma(x_i,x_j)= \delta_{ij}-x_ix_j.\eeq

This operator is invariant under the rotations of $\bbR^{d+1}$, and as a consequence its reversible probability measure is the uniform measure on the sphere (normalized to be a probability). 
A system of coordinates for the upper half sphere $\{x_{d+1}>0\}$ is given by $(x_1, \cdots, x_d)\in \bbB_d$, where $\bbB_d = \{\sum_1^d x_i^2= \|x\|^2< 1\}$ is the unit ball in $\bbR^d$. In this system of coordinates, and thanks to formula~\eqref{eq.density},  one checks easily  that, up to a normalizing constant,
 the reversible  measure is  $(1-\|x\|^2)^{-1/2} dx$, which is  therefore the density of the uniform measure on the sphere in this system of coordinates (see \cite{bglbook}).
 
Now, one may consider some larger dimension $m>d$ and project the Brownian motion on $\bbS^{m}$ on the unit ball in $\bbR^d$ through 
$(x_1, \cdots, x_{m+1})\mapsto (x_1, \cdots, x_d)$. Formula~(\ref{lapl.sphere}) provides immediately that this image is again a diffusion process with generator
\beq\label{sphere.projetee} \cL^{(m)}(x_i)= -mx_i, \Gamma(x_i,x_j)= \delta_{ij}-x_ix_j, \eeq that is the same formula as ~\eqref{lapl.sphere}  except that now $m$ is no longer the dimension of the ball. Once again, formula~(\ref{eq.density}) provides the reversible  measure for this operator, which is, up to a normalizing constant,
$(1-\|x\|^2)^{(m-1-d)/2} dx$, which is therefore the image measure of the uniform measure of the sphere through this projection.

As before, the boundary of the domain (the unit ball)  has equation $\{1-\|x\|^2=0\}$, and we have a boundary equation 
\beq\label{eq.bord.sphere} \Gamma\big(x_i, \log(1-\|x\|^2)\big)= -2x_i,\eeq
similar to equation~(\ref{eq.bord.deltoide}).

Now, it is  again easily checked that, for the function $h= (1-\|x\|^2)^{-(m-1-d)/2}$, one has 
$$\cL^{(m)}(h)=d(m-d-1)h, $$
so that one may perform the associated $h$-transform for which
$$(\cL^{(m)})^{(h)}= \cL^{(2d+2-m)}.$$ In the case where $m=d$, on sees that $\cL^{(d)}$, which is the Laplace operator in this system of coordinates, is transformed into $\cL^{(d+2)}$, which is the projection of the spherical Laplace operator in $\bbS^{d+2}$ onto the unit ball in $\bbR^d$. 

As a consequence, we get
\bprop
A spherical Brownian motion on the unit sphere  $\bbS^{d}\subset \bbR^{d+1}$ conditioned to remain in a half sphere $\{x_{d+1}>0\}$,  has the law of the projection of a spherical Brownian motion on $\bbS^{d+2}$ onto the unit ball in $\bbR^{d}$, lifted on the half upper sphere in $\bbR^{d+1}$. 
\eprop
\section{General $h$-transform for models associated with orthogonal polynomials\label{sec.orthogonal.poly}}

We shall see in this section that all the above examples appear as particular examples, or limit examples, of a very generic one when orthogonal polynomials come into play. Everything relies on a boundary equation similar to~(\ref{eq.bord.deltoide}) or~(\ref{eq.bord.sphere}), which appears as soon as one has a family of orthogonal polynomials which are eigenvectors of diffusion operators.

Let us recall some basic facts about diffusion associated with orthogonal polynomials, following~\cite{BOZ2013}.  We are interested in bounded open sets $\Omega\subset \bbR^d$, with piecewise $\cC^1$ boundary.  On $\Omega$,   we have a probability measure $\mu$ with smooth density $\rho$ with respect to the Lebesgue measure, and  an elliptic  diffusion operator $\cL$ which is symmetric in $\bbL^2(\mu)$. We suppose moreover that polynomials belong to the domain of $\cL$,  and that $\cL$ maps the set $\cP_k$ of polynomials with total  degree less than $k$ into itself.  Then, we may find a $\bbL^2(\mu)$ orthonormal basis formed with polynomials which are eigenvectors for $\cL$.   Following \cite{bglbook}, this is entirely described by the triple $(\Omega, \Gamma, \mu)$, where $\Gamma$ is the square field operator of $\cL$.

We call  such a system $(\Omega, \Gamma, \mu)$ a polynomial system.

Then, one of the main results of~\cite{BOZ2013}  is the following
\bthm\label{thm.polyn.system}~
\benum

\item The boundary $\partial \Omega$ is included in an algebraic surface  with reduced equation $\{P=0\}$, where $P$ is a polynomial which may we written as $P_1\cdots P_k$, where the polynomials $P_i$  are real,  and complex irreducible.

\item If $\cL= \sum_{ij} g^{ij} \partial^2_{ij} + \sum_i b^i \partial_i$, where the coefficients $g^{ij}$ are degree  at most $2$ polynomials and $b^i$ are  polynomials with degree at most $1$. 

\item  The polynomial $P$ divides $\det(g^{ij})$ (that we write $\det(\Gamma)$ in what follows, and which is a polynomial with degree at most $2d$).

\item For each irreducible  polynomial $P_r$ appearing in the equation of the  boundary, there exist  polynomials $L_{i,r}$ with degree at most $1$ such that
\beq\label{eq.g.polyn}\ \forall i= 1, \cdots, d, ~\sum_{j} g^{ij} \partial \log  P_r= L_{i,r} .
\eeq

\item Let  $\Omega$ be  a bounded set, with boundary  described  by a reduced polynomial equation $\{P_1\cdots P_k=0\}$, such that there exist a solution $(g^{ij}, L_{i,k})$ to equation~(\ref{eq.g.polyn}) with $(g^{ij})$ positive definite in $\Omega$. Call $\Gamma(f,f)= \sum_{ij} g^{ij} \partial_i f\partial_j f$ the associated squared field operator. Then for any choice of real numbers $\{\alpha_1, \cdots, \alpha_k\}$ such that $P_1^{\alpha_1}\cdots P_k^{\alpha_k}$ is integrable over $\Omega$ for the Lebesgue measure, setting \\
$$\mu_{\alpha_1, \cdots, \alpha_k}(dx) = C_{\alpha_1, \cdots, \alpha_k} P_1^{\alpha_1}\cdots P_k^{\alpha_k} dx,$$ where $C_{\alpha_1, \cdots, \alpha_k} $ is a normalizing constant, then $(\Omega, \Gamma, \mu_{\alpha_1, \cdots, \alpha_k})$ is a polynomial system.

\item When $P= C\det(\Gamma)$, that is when those 2 polynomials have the same degree, then there are no other measures $\mu$ for which $(\Omega, \Gamma, \mu)$ is a polynomial system.
\eenum 
\ethm 

\brmq\label{rmk.compact.form} Equation~(\ref{eq.g.polyn}), that we shall call the boundary equation (not to be confused with the equation of the boundary), may be written in a more compact form $\Gamma(x_i, \log P_r)= L_{i,r}$. Thanks to the fact that each polynomial $P_r$ is irreducible, this is also equivalent to the fact that $\Gamma(x_i, \log P)= L_i$, for a family $L_i$ of polynomials with degree at most 1.

\ermq

One must be a bit careful about the reduced equation of the boundary $\{P=0\}$, when $P= P_1\cdots P_k$.  This means that each regular point of the boundary is contained in exactly one of the algebraic surfaces $\{P_i(x)=0\}$, and that for each $i= 1\cdots k$, there is at least one regular point  $x$ of the boundary such that $P_i(x)= 0$. In particular,  for a regular point $x\in \partial \Omega$ such that $P_i(x)=0$,  then for $j\neq i$, $P_j(x) \neq 0$ in a neighborhood $\cU$ of such a point, and $P_i(x)=0$ in  $\cU\cap \partial \Omega$.  It is not too hard to see that such a polynomial $P_i$, if real irreducible, is also complex irreducible (if not, it would be written as $P^2+ Q^2$, and $P=Q=0$ on  $\cU\cap  \partial \Omega$).
It is worth to observe that since $P$  divides $\det(\Gamma)$ and that $(g^{ij})$ is positive definite on $\Omega$, then no one of the polynomials $P_i$ appearing in the boundary equation may vanish in $\Omega$. We may therefore chose them to be all  positive on $\Omega$.

The reader should also  be aware that equation~(\ref{eq.g.polyn}), or more precisely the compact form given in Remark~\ref{rmk.compact.form}, and which is 
 the generalization of equations~(\ref{eq.bord.deltoide}) and~(\ref{eq.bord.sphere}),  is a very strong constraint on the polynomial $P$. Indeed,  given $P$, if  one wants to determine the coefficients $(g^{ij})$ and $L_{i}$, this equation is a linear equation in terms of the coefficients of $g^{ij}$ and $L_{i}$, for which we expect to find some non vanishing solution. But the number of equations is much bigger than the number of unknowns, and indeed very few polynomials $P$ may satisfy those constraints. In  dimension $2$ for example, up to affine invariance, there are exactly $10$ such polynomials, plus one one parameter family (see~\cite{BOZ2013}).  The deltoid curve of paragraph~\ref{subsec.deltoid} is just one of them.

\brmq
We shall not use the full strength of this theorem in the  examples developed here. The important fact is the boundary equation~(\ref{eq.g.polyn}), which may be checked directly on many examples, and is the unique property required for the general $h$-transform described in Theorem~\ref{thm.polyn.models}. 
\ermq

Given a bounded set $\Omega$ and an operator $\Gamma$  satisfying the conditions  of Theorem~\ref{thm.polyn.system},  and for any choice of  $\{\alpha_1, \cdots, \alpha_k\}$ such that $P_1^{\alpha_1}\cdots P_k^{\alpha_k}$ is integrable over $\Omega$ for the Lebesgue measure, we have a corresponding symmetric operator  $\cL_{\alpha_1, \cdots, \alpha_k}$.  For this operator, as was the case in paragraphs~\ref{subsec.deltoid} and~\ref{subsec.spheres}, one may extend the integration by parts~(\ref{eq.ipp}) to any pair of polynomials, and this provides a sequence of orthogonal polynomials  which are eigenvectors of the operator $\cL_{\alpha_1, \cdots, \alpha_k}$. 

Conversely, the boundary equation~(\ref{eq.g.polyn}) is automatic as soon as we have a generator on a bounded set with regular boundary,  and a complete system  of eigenvectors which are polynomials. But it may happen that those conditions are  satisfied even on non bounded domains, and even when the associated measure is infinite (this appears in general in limits of such polynomial models, as in the Laguerre and Bessel cases). We may therefore give a statement in a quite general setting. 

\bthm\label{thm.polyn.models} Assume that a symmetric positive definite matrix $(g^{ij})$ on some open set $\Omega\subset \bbR^d$, is  such that for any $(i,j)$, $g^{ij}$  is a polynomial of degree at most $2$.  Let us call $\Gamma$ the associated square field operator. Suppose moreover that we have some polynomials $P_k$, positive on $\Omega$, such that, for any $k$, 
\beq\label{eq.bdry.gal}\forall i= 1, \cdots, d,~ \sum_i g^{ij} \partial_j \log P_r= \sum_i \Gamma(x^i, \log P_k) = L_{i,k},\eeq where $L_{i,k}$ are degree $1$ polynomials.  For any $(\alpha_1, \cdots, \alpha_k)$, let $\mu_{\alpha_1, \cdots, \alpha_k}$ be the measure with density $P_1^{\alpha_1}\cdots P_k^{\alpha_k}$ with respect to the Lebesgue measure on $\Omega$, and let 
$\cL_{\alpha_1, \cdots, \alpha_k}$ be the generator associated with the Markov triple $(\Omega, \Gamma, \mu_{\alpha_1, \cdots, \alpha_k})$. 

Then,  there exist constants $c_k$ such that, for any $(\alpha_1, \cdots, \alpha_k)$, the  function $h= P_1^{-\alpha_1}\cdots P_k^{-\alpha_k}$ satisfies
$$\cL_{\alpha_1, \cdots, \alpha_k}(h)=-(\sum_k \alpha_k c_k)h.$$ 
Moreover, 
$(\cL_{\alpha_1, \cdots, \alpha_k})^{(h)}= \cL_{-\alpha_1, \cdots, -\alpha_k}$.

\ethm

\bpf   
We shall prove the assertion with  
$c_k = \sum_i \partial_i L_{i,k}$.

With $\rho= P_1^{\alpha_1}\cdots  P_k^{\alpha_k}$,  we write our operator  $\cL_{\alpha_1, \cdots, \alpha_k}$ as 
$$\sum_{ij} g^{ij} \partial^2_{ij} + \sum_i b^i\partial_i,$$ where 
\beq \label{eq.bi} b_i= \sum_j \partial_j g^{ij} + \sum_{r,j}\alpha_r g^{ij} \partial_j \log P_r= \sum_j \partial_j g^{ij} + \sum_r\alpha _r L_{i,r}.\eeq
With
$$\cL_0= \sum_{ij} g^{ij} \partial^2_{ij} + \sum_i \partial_{j}g^{ij}\partial_i,$$
then
\beq\label{eq.defL0}\cL_{\alpha_1, \cdots, \alpha_k}(f)= \cL_0(f) + \sum_i \alpha_i \Gamma(\log P_i, f).\eeq

What we want to show is 
$\cL_{\alpha_1, \cdots, \alpha_k}(h)= ch$, or 
$$\cL_{\alpha_1, \cdots, \alpha_k}(\log h) +\Gamma(\log h, \log h) = c.$$
With, $\log h= -\sum_i \alpha_i \log P_i$, and comparing with equation~\eqref{eq.defL0}, this amounts to
$$\cL_0(\log h)=-\sum_r \alpha_q \cL_0(\log P_r) =c.$$

We may first  take derivative in equation~(\ref{eq.bdry.gal}) with respect to $x_i$ and add the results in $i$ to get
$$\sum_{ij} g^{ij} \partial_{ij}\log P_r + \sum_i  \partial_i (g^{ij} )\partial_j \log P_r =\sum_i \partial_i L_{i,r}= c_r, $$
that is $\cL_0(\log P_r)= c_r$.

It remains to add these identities  over $r$ to get the required result.

Comparing the reversible  measures, it is then immediate to check that $(\cL_{\alpha_1, \cdots, \alpha_k})^{(h)}= \cL_{-\alpha_1, \cdots, -\alpha_k}$, 
\qed
\epf

\brmq
The function $h$ is always the inverse of the density  with respect to the Lebesgue measure, in the system of coordinates in which we have this polynomial structure. Of course, the choice of the coordinate system is related to the fact that,  in those coordinates,  we have orthogonal polynomials (at least when the measure is finite on a bounded set).  In the Bessel case, for example, which is a limit of a Laguerre models, one has to change $x$ to $x^2$ to get a simple correspondance between the $h$ function and the density. The same is true in many natural examples, where one has to perform some change of variable to get the right representation (for example from the triangle to the deltoid in paragraph~\ref{subsec.deltoid}).

\ermq

\brmq In many situations, there are natural geometric interpretations for these polynomial models when the parameters $(\alpha_1, \cdots , \alpha_k)$ are half integers, in general  with $\alpha_i \geq -1/2$. The case $\alpha_i = -1/2$ often corresponds to Laplace operators, while the dual case  $\alpha_i= 1/2$  often corresponds to the projection of a Laplace operator in larger dimension.

\ermq

\section{Further examples\label{sec.further.examples}}

We shall provide two more examples, one  which follows directly from Theorem~\ref{thm.polyn.models}, and another one on a  non bounded domain with infinite measure. One may provide a lot of such examples, many of them arising from Lie group theory, Dunkl operators, random matrices, etc.  However, we chose to present those two cases because they put forward some specific features of diffusion operators associated with orthogonal polynomials.

\subsection{Matrix Jacobi processes\label{subsec.matrix.jacobi}}

This model had been introduced by Y. Doumerc in his thesis~\cite{YanDoumerc},  and  had also been studied in the complex case, especially from the asymptotic point of view in~\cite{DemniHamdiFreeJacobi12, DemniHmidiFreeUnitaryBM12}. It plays a similar rôle than the one-dimensional Jacobi processes for matrices. One starts from the Brownian motion on the group $SO(d)$.  Since $SO(d)$ is a semi-simple compact Lie group, it has a canonical Casimir operator similar to the one described in equation~\eqref{casimir.SU(d)}. If $O= (m_{ij})$ is an $SO(d)$ matrix, then the  Casimir operator  may be described through it's action on the entries $m_{ij}$. One gets
\beq\label{casimir.S0(d)} \cL(m_{ij})=-(d-1) m_{ij}, \quad 
\Gamma( m_{kl},m_{qp})= \delta_{(kl)(qp)}- m_{kp}m_{ql}.
\eeq
Observe that when restricted to a single line or column, one recovers the spherical Laplace operator on $\bbS^{d-1}$ described in equation~\eqref{lapl.sphere}. 

An $SO(d)$-Brownian matrix is then a diffusion process with generator this Casimir operator on $SO(d)$. 

It is again clear from the form of the operator $\cL$ that it preserves for each $k\in \bbN$ the set of polynomials in the entries $(m_{ij})$ with total degree less that $k$. However, these "coordinates" $(m_{ij})$ are not independent, since they satisfy algebraic relations, encoded in the fact that $OO^*= \Id$. We may not apply directly our main result Theorem~\ref{thm.polyn.models}. We shall nevertheless look at some projected models on which the method applies.

One may extract  some $p\times q$ submatrix $N$  by selecting $p$ lines and $q$ columns, and we observe that the generator acting on  the entries of this extracted matrix  $N$ depend only on the entries of $N$.  Therefore, the operator projects on these extracted $p\times q$ matrices and the associated process is again a diffusion process : we call  this the projection of the Brownian motion in $SO(d)$ onto the set $\cM_{p,q}$ of $p\times q$ matrices. Thanks to formula~\eqref{eq.density}, one may compute the density of the image measure, with respect to the Lebesgue measure in the entries of $N$. Whenever $p+q\leq d$,  it happens to be, up to a normalizing constant
$\det(\Id-NN^*)^{(d-1-p-q)/2}$, with support the set $\Omega= \{N, NN^*\leq \Id\}$.  This formula is easy to check if we recall that, for a matrix $M$ with entries $(m_{ij})$,
$$\partial_{m_{ij}} \log \det(M)= M^{-1}_{ji},$$ a consequence of Cramer's formula.

When $p+q\geq d+1$,  there are however algebraic relations between the entries of  $N$ and the image measure has no density with respect to the Lebesgue measure. 
For example, when $p+q= d+1$, then the measure concentrates on the algebraic set $\{ \det(\Id- NN^*)= 0\}$.  It may be checked that it has a density with respect of the Lebesgue measure of this hypersurface.
Indeed, one may fix $p$ and $q$ and consider $d$  as a parameter. It is worth to observe that the function $\det(\Id-NN^*)^\alpha$ is not integrable on the domain $\Omega$ whenever $\alpha \leq -1$. Moreover, , when  $\alpha>-1$ and $\alpha\to -1$,  the probability measure  with density $ C_\alpha \det(\Id-NN^*)^{\alpha}$ concentrates on the set $\{\det(\Id-NN^*)= 0\}$, and the limit is a measure supported by this surface  with a density with respect of the  surface  measure. 
Things become even worse as the number $p+q$ increases, the measure being concentrated on manifolds with higher and higher co-dimensions. 

We are in a situation different from the sphere case here, since we may not chose  the parameters in which the operator has a nice polynomial expression as a local system of coordinates.  Indeed, 
the Lie group $SO(d)$ is a $d(d-1)/2$ manifold.  Since we want  algebraically independent coordinates,   we are limited to $pq$ ones,  with $p+q\leq d$, we may  have at most $d^2/4$ algebraically independent such polynomial coordinates, which for  $d>2$   is less than the dimension of the manifold.

It is worth to observe that, again when $p+q\leq d$, one has $pq$ variables, the determinant of the metric $\Gamma$ is a degree $2pq$ polynomial, whereas $\det(\Id-NN^*)= \det(\Id-N^*N)$ is of degree at most $2\min(p,q)$. We are not in  the case of maximal degree for the boundary equation. When $p+q= d$, the density measure is $\det(\Id-NN^*)^{-1/2}$, but the corresponding operator is not a Laplace operator (for which the density of the  measure would be $\det(\Gamma)^{-1/2}$).  Since we are in the situation of orthogonal polynomials as described in Section~\ref{sec.orthogonal.poly}, we know that we may perform an $h$-transform.  

For the particular case where  $d= p+q$, we get
\bprop The matrix $N$ projected  from an $SO(d)$-Brownian  matrix on $\cM_{p,q}$ conditioned to remain in the set $\{NN^*<  \Id\}$ has  the  law of the projection of  a $SO(d+2)$-Brownian  matrix on $\cM_{p,q}$.
\eprop

\subsection{ Brownian motion in a Weyl chamber\label{sec.Weyl.chamber}}

This last example is again quite well known, but is happens to fit also with the general picture associated with orthogonal polynomials, although no orthogonal polynomials are associated with it.  Indeed, it does  not follow directly from the  setting of  Section~\ref{sec.orthogonal.poly}, one the one side because it is non compact, on the other because the reversible measure in this situation is infinite.   But it satisfies the all the algebraic properties described is Section~\ref{sec.orthogonal.poly}, and we may then check that we may apply   the result for the associated $h$-transforms. Indeed, one may replace in what follows Brownian motion by Ornstein-Uhlenbeck operators, which have as reversible  measure a Gaussian measure with variance $\sigma^2$, and then let $\sigma$ go to infinity. In the Ornstein-Uhlenbeck case, we are in the setting  of orthogonal polynomials,   however with a non bounded domain. But this would introduce further complication, since the Brownian case gives simpler formulas.

As described above, the $h$-transform is easy to compute in a system of coordinates which have some relevant polynomial structure. Here, one  good choice for the coordinate system are the elementary symmetric functions in $d$ variables. We shall perform mainly computations on these elementary symmetric functions of the components of the $d$-dimensional Brownian motion, following~\cite{BakryZaniClifford}.  In $\bbR^d$, one may consider the Brownian $(B^1_t, \cdots, B^d_t)$ and reflect it around the hyperplanes which are the boundaries of the set  $\{x_1< \cdots < x_d\}$, which is usually called a Weyl chamber. To describe this reflected Brownian motion, it is easier to consider the elementary symmetric functions which are the coefficients of the polynomial 
$$P(X)= \prod_{i=1}^d (X-x_i)= \sum_{i= 0}^d a_i X^i,$$ where 
$a_d=1$ and the functions $a_i, i= 0, \cdots d-1$ are, up to a sign,  the elementary symmetric functions of the variables $(x_i)$. The map $(x_i)\mapsto (a_i)$ is a diffeomorphism in the Weyl chamber $\{x_1 < \cdots< x_d\}$ onto it's image. To understand the image, one has to consider the discriminant $\discrim(P)$,  a polynomial in the variables $(a_i)$,  which is, up to a sign $(-1)^{d(d-1)/2}$,  the following $(2d-1)\times (2d-1)$ determinant

\beqnas\bpm 	1& a_{d-1} & a_{d-2}&\cdots & a_0&0&\cdots&0\\
		0&1         	& a_{d-1}	&\cdots & a_1&a_0&\cdots &0\\
		0&0		&1        	 & \cdots     &a_2&a_1&\cdots & 0\\
		\cdots & \cdots &\cdots                   &\cdots           & \cdots      &  \cdots     & \cdots&\cdots\\
		\cdots & \cdots &     \cdots              &      \cdots     &    a_{p-2}   &  \cdots     & a_1&a_0\\
		1& (d-1)a_{d-1} & (d-2)a_{d-2}&\cdots & a_1&0&\cdots&0\\
		0&1         	& (d-1)a_{d-1}	&\cdots & 2a_2&a_1&\cdots &0\\
		0&0		&1        	 & \cdots     &3a_3&2a_2&\cdots & 0\\
		\cdots & \cdots &\cdots                   &\cdots           & \cdots      &  \cdots     & \cdots&\cdots\\
		\cdots & \cdots  &      \cdots              &      \cdots      &     \cdots   &  \cdots     & 2a_2&a_1
\epm
\eeqnas

It turns out that this discriminant is  $\prod_{i<j}(x_j-x_i)^2$. The image of the Weyl chamber is the connected component $\Omega$ of the set $\{\discrim(P) \neq 0\}$ which contains the image of any polynomial with $d$ real distinct roots, and the image of the boundary of the Weyl chamber is  $\partial \Omega$, a subset of  the algebraic surface $\{\discrim(P)=0\}$.  It is not hard to observe (by induction on the dimension $d$) that the image of the Lebesgue measure $dx_1\cdots dx_d$ on the Weyl chamber is nothing else than $\un_\Omega\discrim(P)^{-1/2} \prod da_i$.  

Now, the Brownian motion in $\bbR^d$ may be described, up to a factor 2,  through 
$$\Gamma(x_i, x_j)= \delta_{ij}, ~\Delta(x_i)=0.$$ We want to describe this operator acting on the variables $(a_0, \cdots, a_{d-1})$. 
Since any of the functions $a_j$ is a degree $1$ polynomial in the variables $x_i$, one has $\Delta(a_j)=0$,  $j= 0, \cdots, d$. To compute $\Gamma(a_i,a_j)$, it is simpler to compute 
$$\Gamma(P(X), P(Y))= \sum_{i,j} X^i Y^j \Gamma(a_i,a_j).$$ We obtain
\bprop The image of the operator $\Delta$ in $\bbR^n$ on the coefficients of the polynomial  $P(X)= \prod_i (X-x_i)$ is given by
\beq \label{Gamma.Weyl.chamber}\Gamma(P(X), P(Y))= \frac{1}{Y-X}\big(P'(X)P(Y)-P'(Y)P(X)\big), ~\Delta(P(X))=0.\eeq
\eprop
\bpf The second formula is a direct consequence of $\Delta(a_i)=0$, while for the first, it is simpler to look at $\Gamma(\log P(X), \log P(Y))$. 
\beqnas \Gamma(\log P(X), \log P(Y))&=& \sum_{ij} \Gamma(\log (X-x_i), \log (Y-x_j))\\&=& \sum_{ij} \frac{1}{(X-x_i)(Y-x_j)} \Gamma(x_i,x_j)
\\&=& \sum_{i}  \frac{1}{(X-x_i)(Y-x_i)}= \frac{1}{Y-X}\Big(\frac{P'(X)}{P(X)}-\frac{P'(Y)}{P(Y)}\Big).
\eeqnas 
\qed
\epf 

\brmq From formula~\eqref{Gamma.Weyl.chamber}, it is clear that $\Gamma(a_i,a_j)$ are degree $2$ polynomials in the variables $a_i$.
\ermq

The image of the Brownian motion $B_t$ in the variables $(a_i)$ is nothing else than the Brownian motion reflected through the walls of the Weyl chamber. Its generator is described through the $\Gamma$ operator given in equation~(\ref{Gamma.Weyl.chamber})  and it is the image of the Laplace operator on the Weyl chamber. Since it is an Euclidean Laplace operator,   the reversible measure is, up to a constant, $\det(\Gamma)^{-1/2}$,  and this shows that the determinant $\det(\Gamma)$ of the  metric is, up to a constant, $\discrim(P)$. 

Moreover, from the general representation of diffusion operators ~(\ref{diff.coord}), and the equation~(\ref{eq.density}) giving the reversible  measure, we have,  with $\rho= \discrim(P)^{-1/2}$, $b_i=0$,
\beq\label{eq.bdry.weyl}\sum_{ij} \Gamma(a_i,a_j)\partial_{a_i} \log \rho = -\sum_j \partial_{a_j} \Gamma(a_i,a_j).\eeq
Since $ \partial_{a_j} \Gamma(a_i,a_j)$ is a degree at most one polynomial in the variables $a_i$, this is nothing else than the boundary equation~(\ref{eq.g.polyn}) for general polynomial models. We may therefore apply the general result described in Section~\ref{sec.orthogonal.poly}.

In order to identify the result of the $h$-transform, an important formula relating $\Gamma$ and the discriminant function is the following
\bprop\label{prop.eq.discrim.Gamma} For the operator $\Gamma$ defined in~\eqref{Gamma.Weyl.chamber}, one gas\beq\label{eq.discrim.P}\Gamma\big(P(X), \log \discrim(P)\big)= -P''(X).\eeq
\eprop

\bpf
One may find a proof of this formula in~\cite{BakryZaniClifford}, but the one we propose here is simpler. To check equation~(\ref{eq.discrim.P}), it is enough to establish it it in a Weyl chamber $\{x_1<x_1<\cdots< x_d\}$ where $P(X)= \prod(X-x_i)$ and $\discrim(P)= \prod_{i<j} (x_i-x_j)^2$, since the map $(x_1, \cdots, x_k)\mapsto P(X)$ is a local diffeormorphism in this domain.

In those coordinates, $\Gamma(x_i, x_j)= \delta_{ij}$
and, from the change of variable formula~(\ref{change.var.Gamma}), one has 
$$\Gamma(\log P(X), \discrim (P))= 2\sum_{i, j<k} \Gamma(\log (X-x_i), \log (x_j-x_k))= -2\sum_{i,j<k}\frac{1}{X-x_i}\frac{1}{x_j-x_k}(\delta_{ij}-\delta_{ik}).$$
From which one gets
$$\Gamma(\log P(X), \log \discrim(P))=-2\sum_{i\neq j} \frac{1}{X-x_i}\frac{1}{x_i-x_j}.$$
On the other hand,
\beqnas \frac{P''}{P}&=& \Big( \frac{P'}{P}\Big)'+ \Big(\frac{P'}{P}\Big)^2= \sum_{i\neq j} \frac{1}{(X-x_i)(X-x_j)}\\
&=& \sum_{i\neq j}\big( \frac{1}{X-x_i}-\frac{1}{X-x_j}\big)\frac{1}{x_i-x_j},= 2\sum_{i\neq j} \frac{1}{X-x_i}\frac{1}{x_i-x_j}.\eeqnas
From this we get
$$\Gamma(\log P, \log \discrim(P))= -\frac{P''}{P},$$ which in turns gives~(\ref{eq.discrim.P}).
\qed\epf

Proposition~\ref{prop.eq.discrim.Gamma} is central in the identification of various processes with the same $\Gamma$ given by~\eqref{Gamma.Weyl.chamber}. It turns out that the same operator with this  $\Gamma$ operator and reversible  measure $\discrim(P)^{1/2}$ has a nice geometric interpretation: namely, it is the Dyson complex  process, that is the law of the spectrum of Hermitian Brownian matrices, introduced by Dyson~\cite{dyson}. In the same way, the case where the reversible  measure is the  Lebesgue measure corresponds to Dyson process for real symmetric matrices, and $\rho= \discrim(P)^{3/2}$ corresponds to Dyson process for symmetric quaternionic matrices~, see~\cite{AndGuioZeit, BakryZaniClifford, Forrester}.

Let us show a direct way to check this (first in the real symmetric case, where it is simpler).
The Brownian motion on symmetric matrices is nothing else that the Brownian motion of the Euclidean space of symmetric matrices $M$, endowed with the  Euclidean norm $\| M\|^2= \tr(M^2)$. When $M= (m_{ij})$, this may be described as 

 $$\Gamma(m_{ij}, m_{kl})= \frac{1}{2} (\delta_{ik}\delta_{jl}+ \delta_{il}\delta_{jk}), ~{\cL}(m_{ij})=0.$$ 
 One may look at its action of the characteristic polynomial $P(X)= \det(X\Id-M)$.  We get
 \bprop\label{prop.spect.real.symm} For the characteristic polynomial associated with a Brownian symmetric  matrix, one has
 $$
 \Gamma\big(\log P(X), \log P(Y)\big)= \frac{1}{Y-X}\Big(\frac{P'(X)}{P(X)}-\frac{P'(Y)}{P(Y)}\Big), 
~ \cL P(X)= -\frac{1}{2}P''.$$
 \eprop
 \bpf
 To compute $\Gamma\big(P(X), P(Y)\big)$ and $\cL \big(P(X)\big)$. In order to apply the change of variable formula~(\ref{diff.coord}), we may apply the general formulas for the determinant function
 $$\partial_{m_{ij}} \log \det~ M= M^{-1}_{ji}, ~\partial_{m_{ij}}\partial_{m_{kl}}\log \det~ M= -M^{-1}_{jk}M^{-1}_{li},$$ which are direct consequences of Cramer's formulas for the inverse matrix.
 
 Then the formulas are direct applications of the chain rule formula.\qed
 
\epf

We may now compare with equation~(\ref{eq.discrim.P}) to see that the reversible  measure for the spectral measure for Brownian symmetric matrices, given by the general formula~(\ref{eq.density}), in the system of coordinates which are the coefficients  $(a_i)$ of the characteristic polynomial, is  the Lebesgue measure.

We may perform the same computation for Hermitian matrices. In this situation, one would consider a complex valued matrix $M$ with entries $(z_{ij})$ and satisfying
$$\Gamma(z_{ij}, z_{kl})= 0,~ \Gamma(z_{ij},\bar  z_{kl})=\delta_{il}\delta_{jk},  ~\cL(z_{ij})=0.$$

One may again perform the same computation on  $P(X)= \det(X\Id-M)$, and we get
\bprop For the characteristic polynomial associated with a Brownian Hermitian   matrix, one has
 $$
 \Gamma(P(X), P(Y))= \frac{1}{Y-X}\big( P'(X)P(Y)-P'(Y)P(X)\big), 
~ \cL P(X)= -P''.$$
 \eprop

We do not give the proof, which follows along the same lines that the one of Proposition~\ref{prop.spect.real.symm}. More details may be found in~\cite{BakryZaniClifford}.

  As a consequence, comparing with equation~\eqref{eq.discrim.P} and  equation~\eqref{eq.density} in the system of coordinates given by the coefficients of $P(X)$, the density of the reversible measure for the Hermitian Dyson process is $\discrim(P)^{1/2}$ whereas  the density of the  reversible measure of the Brownian motion in the Weyl chamber is $\discrim(P)^{-1/2}$.

 Transfering back to the Weyl Chamber through the local diffeomorphism between the coefficients of $P(X)$ and the roots $(x_1<x_2<\cdots< x_d)$ of $P(X)$. We  obtain
\bprop
 The Brownian motion conditioned not to reach the boundary of the Weyl chamber $\{x_1<\cdots< x_d\}$ has the law of the spectrum of an Hermitian $d\times d$ matrix.
\eprop

 \bibliographystyle{amsplain}   
\bibliography{bib.h.tranforms}

\end{document}